\documentclass{amsart}
\usepackage{amscd,amssymb,enumerate,amsfonts,amsmath,verbatim}

\newcommand{\id}{\operatorname{id}}

\newcommand{\codim}{\operatorname{codim}}
\newcommand{\codepth}{\operatorname{codepth}}
\newcommand{\edim}{\operatorname{edim}}
\newcommand{\depth}{\operatorname{depth}}
\newcommand{\card}{\operatorname{card}}

\newcommand{\rank}{\operatorname{rank}}

\newcommand{\shift}{{\Sigma}}
\newcommand{\filt}{\operatorname{F}}

\newcommand{\Ext}{\operatorname{Ext}}
\newcommand{\Tor}{\operatorname{Tor}}
\newcommand{\Hom}{\operatorname{Hom}}
\newcommand{\HH}{\operatorname{H}}
\newcommand{\Ker}{\operatorname{Ker}}
\newcommand{\Coker}{\operatorname{Coker}}
\newcommand{\Imag}{\operatorname{Im}}
\newcommand{\Po}{{P}}
\newcommand{\Ba}{{I}}

\newcommand{\xra}{\xrightarrow}
\newcommand{\lra}{\longrightarrow}

\newcommand{\Ho}{{\mathbf{H}}}
\newcommand{\B}{{\mathbf{B}}}
\newcommand{\TE}{{\mathbf{TE}}}
\newcommand{\G}{{\mathbf{G}}}
\newcommand{\GO}{{\mathbf{GO}}}

\newcommand{\Spec}{\operatorname{Spec}}

\newcommand{\Ass}{\operatorname{Ass}}

\newcommand{\Supp}{\operatorname{Supp}}
\newcommand{\height}{\operatorname{height}}

\newcommand{\lf}{\otimes^{\mathbf L}}
\newcommand{\rh}{{\mathbf R}\!\Hom}

\newcommand{\dd}{\partial}

\newcommand{\BZ}{{\mathbb Z}}

\newcommand{\fa}{{\mathfrak a}}
\newcommand{\fb}{{\mathfrak b}}
\newcommand{\fc}{{\mathfrak c}}
\newcommand{\fd}{{\mathfrak d}}

\newcommand{\fm}{{\mathfrak m}}
\newcommand{\fn}{{\mathfrak n}}
\newcommand{\fp}{{\mathfrak p}}
\newcommand{\fq}{{\mathfrak q}}
\newcommand{\fr}{{\mathfrak r}}

\newcommand{\Sym}{{\mathsf S}}

\newcommand{\EE}[3]{{}^{#1\!}\operatorname{E}_{#2,#3}}

\newcommand{\DD}[3]{{}^{#1\!}d_{#2,#3}}
\newcommand{\du}[2]{D_{#1}(#2)}

\newcommand{\dua}[1]{D^{#1}}

\newcommand{\wh}{\widehat}

\newcommand{\col}{\colon}

\newcommand{\bsf}{{\boldsymbol f}}
\newcommand{\bsg}{{\boldsymbol g}}

\newcommand{\bsx}{{\boldsymbol x}}

\newcommand{\ges}{\geqslant}
\newcommand{\les}{\leqslant}

\theoremstyle{remark}
\newtheorem{step}{Step}

\swapnumbers

\theoremstyle{plain}
\newtheorem{theorem}{Theorem}[section]

\newtheorem{proposition}[theorem]{Proposition}
\newtheorem{lemma}[theorem]{Lemma}
\newtheorem{corollary}[theorem]{Corollary}

\theoremstyle{definition}
\newtheorem{chunk}[theorem]{}
\newtheorem{subchunk}{}
\newtheorem*{chunk*}{}

\theoremstyle{remark}

\newtheorem*{Remark}{Remark}
\newtheorem{remark}[theorem]{Remark}
\newtheorem*{Conjecture}{Conjecture}
\newtheorem*{Question}{Question}

\numberwithin{equation}{theorem}
\numberwithin{subchunk}{theorem}

\begin{document}
\title[Extensions of a dualizing complex]{Extensions of a dualizing
complex by its ring:\\ Commutative versions of\\a conjecture of
Tachikawa} \author[L.~L.~Avramov]{Luchezar L.~ Avramov}
\address{Department of Mathematics, University of Nebraska, Lincoln,
NE 68588, U.~S.~A.}
\email{avramov@math.unl.edu}
\author[R.-O.~Buchweitz]{Ragnar-Olaf Buchweitz}
\address{Department of Mathematics, University of Toronto, Toronto
ON M5S 3G3, Canada}
\email{ragnar@math.utoronto.ca}
\author[L.~M.~\c Sega]{Liana M.~\c Sega}
\address{Department of Mathematics, Purdue University, West
Lafayette, In 47907, U.~S.~A.}
\curraddr{Mathematical Sciences Research Institute, 1000 Centennial Drive, Berkeley, CA 94720}
\email{lmsega@math.purdue.edu}

 \begin{abstract}
Let $(R,\fm,k)$ be a commutative noetherian local ring with dualizing
complex $\dua R$, normalized by $\Ext^{\depth(R)}_R(k,\dua R)\cong k$.
Partly motivated by a long standing conjecture of Tachikawa
on (not necessarily commutative) $k$-algebras of finite
rank, we conjecture that if $\Ext^n_R(\dua R,R)=0$ for all
$n>0$, then $R$ is Gorenstein, and prove this in
several significant cases.
  \end{abstract}

\dedicatory{Dedicated to Wolmer Vasconcelos on the occasion of his
65th birthday}

\thanks{L.L.A.\ was partly supported by a grant from the NSF.\endgraf
R.O.B.\ was partly supported by a grant from NSERC.\endgraf
L.M.\c S.\ was a CMI Liftoff Mathematician for 2002.  She thanks
the University of Nebraska-Lincoln for hospitality during the Spring
semester of 2002, and the University of Toronto for hosting a visit.}

\date{\today}
 \maketitle

\section*{Introduction}

Let $(R,\fm,k)$ be a {\em local ring\/}, that is, a commutative noetherian
ring $R$ with unique maximal ideal $\fm$ and residue field $k=R/\fm$.
We write $\edim R$ for the minimal number of generators of $\fm$ and
set $\codepth R=\edim R-\depth R$.

 We assume that $R$ has a dualizing complex, and let $\dua{R}$ denote
such a complex shifted so that $\HH_i(\dua{R})=0$ for $i<0$ and
$\HH_0(\dua{R})\ne 0$; see Section \ref{Dualizing complexes} for
more details on dualizing complexes.  The ring $R$ is Cohen-Macaulay
(respectively, Gorenstein) if and only if $\dua R$ can be taken to be
a dualizing module for $R$ (respectively, to be the $R$-module $R$).
Thus, when $R$ is Gorenstein, $\Ext^i_R(\dua R,R)=0$ holds trivially
for all $i>0$.  We prove that vanishing of such $\Ext$ groups often
implies that $R$ is Gorenstein, as follows:
 \begin{enumerate}[\rm\text{Section} 1.]
  \item[\rm Section \ref{Generically}.]
$R$ is generically Gorenstein (e.g., $R$ is reduced),
and $\Ext_R^i(\dua R,R)=0$ for all $i\in[1,\dim R]$.
  \item[\rm Section \ref{Residue}.]
$R$ is a homomorphic image of a generically Gorenstein ring with
dualizing complex, and $\Ext_R^i(\dua R,R)=0$ for all $i\in[1,\dim
R+1]$.
  \item[\rm Section \ref{Well}.]
$R$ is Cohen-Macaulay, in the linkage class of a complete intersection
(e.g., $\codepth R\le 2$), and $\Ext^i_R(\dua R,R)=0$ for all $i\in[1,\dim
R+1]$.
  \item[\rm Section \ref{Loewy}.]
$R$ contains a regular sequence $\bsf$ with the property that $\fm^3\subseteq(\bsf)$, and
$\Ext_R^i(\dua R,R)=0$ for all $i\in[1,\dim R+1]$.
  \item[\rm Section \ref{Artinian}.]
$R$ is the special fiber of a finite flat local homomorphism that possesses
a Gorenstein fiber, and $\Ext_R^1(\dua R,R)=0$.
  \item[\rm Section \ref{Codepth}.]
$\codepth R\le 3$ and $\Ext^i_R(\dua R,R)=0$ for all $i\gg 0$.
  \item[\rm Section \ref{Golod}.]
$R$ is Golod and $\Ext^i_R(\dua R,R)=0$ for all $i\gg 0$.
 \end{enumerate}

An effort to make sense of the disparate hypotheses leads us to pose the

 \begin{Question}
Does the vanishing of $\Ext^i_R(\dua R,R)=0$ for $(\dim R+1)$
consecutive positive values of $i$ imply that $R$ is Gorenstein?
  \end{Question}

A result of Foxby \cite{Foxby-notes}, see Section \ref{Dualizing
complexes}, implies that if $\Ext_R^i(\dua R,R)$ vanishes for
$i=1,\dots,\dim R$, then $R$ is Cohen-Macaulay, so this hypothesis is
implicit in the

 \begin{Conjecture}
If $\Ext^i_R(\dua R,R)=0$ for all $i>0$,
then $R$ is Gorenstein.
 \end{Conjecture}

This is an analog of a 30 years old conjecture of Tachikawa \cite{Ta} on
(possibly non-commutative) algebras of finite rank over a field.  In fact,
if the hypotheses of both conjectures apply, then they are equivalent.
The result of Section \ref{Loewy} is obtained by transposing the proof,
due to Asashiba and Hoshino \cite{As}, \cite{AsHo}, of one case of
Tachikawa's conjecture.  On the other hand, our results prove that
conjecture in new cases.  Background material and more details on
Tachikawa's conjecture are given in Section \ref{Relations}.

Our arguments use a substantial amount of homological algebra for
complexes of modules.  Appendix \ref{Constructions} contains a synopsis
of the basic constructions.  In Appendix \ref{Lemmas} we have collected
miscellaneous results on complexes that play central roles in various
proofs, but do not draw upon the framework of noetherian local algebra.

\section{Dualizing complexes}
\label{Dualizing complexes}

Let $(R,\fm, k)$ be a local ring. In this section we collect basic
properties of dualizing complexes. A comprehensive treatment can
be found in Foxby's notes \cite{Foxby-notes}.  We refer to Appendix
\ref{Constructions} for basic notation regarding complexes.

A complex of $R$-modules is said to be {\em dualizing\/} if it has finite
homology and there is an integer $d$ such that $\Ext_R^d(k,D)\cong k$
and $\Ext^i_R(k,D)=0$ for $i\ne d$.
 For completeness, we recall
that this is the case if and only if $R$ is a homomorphic image
of a Gorenstein ring, where the ``if'' part is a deep recent result of
Kawasaki \cite{Ka}.  In this paper, we always assume
that $R$ has a dualizing complex.

\begin{chunk}
Let $\dua R$ denote a dualizing complex with $\inf\HH(\dua{R})=0$, and set
$d=\depth R$.

\begin{subchunk} Every dualizing complex is quasi-isomorpic to a shift
of $\dua R$, cf.\ \cite[15.14]{Foxby-notes}.
\end{subchunk}

\begin{subchunk}
\label{size of dc}
$\sup\HH(\dua{R})=\dim R-d$, cf.\ \cite[15.18]{Foxby-notes}.
\end{subchunk}

\begin{subchunk}
\label{id of dc}
$\Ext_R^i(k,\dua R)=0$ for $i\ne d$ and $\Ext_R^d(k,\dua R)\cong k$,
cf.\ \cite[15.18]{Foxby-notes}.
 \end{subchunk}

\begin{subchunk}
\label{localdc}
For every $\fp\in\Spec R$ and for $n=\depth R-\depth R_\fp-\dim (R/\fp)$
there is a quasi-isomorphism $\shift^n(\dua{R})_{\fp}\simeq \dua{R_\fp}$ of
complexes of $R_\fp$-modules, cf.\ \cite[15.17]{Foxby-notes}.
 \end{subchunk}
\end{chunk}

The next result slightly extends an observation by the referee.

\begin{proposition}
\label{associated}
If $\fp$ is an associated prime ideal of $R$, then
\[
\Ext^i_R(\dua R,R)\ne 0 \quad \text{for} \quad i=\dim (R/\fp)-\depth R
\]
\end{proposition}

\begin{proof} Since $D^R$ has finite homology, $\Ext_R(D^R,R)$
localizes, cf.\ \cite[6.47]{Foxby-notes}. This justifies the first
isomorphism below:
\begin{align*}
\Ext^i_R(\dua R,R)_\fp&\cong\Ext^i_{R_\fp}\big((\dua
{R})_\fp,R_\fp\big)\cong\Ext^i_{R_\fp}\big(\shift^{i}(\dua
{R_\fp}),R_\fp\big)\\
&\cong \Ext_{R_\fp}^0(\dua{R_\fp}, R_\fp)
\end{align*}
The second isomorphism comes from \ref{localdc} and the last one is
due to the shift. Now $\inf
\HH(\dua{R_\fp})=0$ by \ref{size of dc}, and hence it follows from
\cite[6.42(3)]{Foxby-notes} that
\[
 \Ext_{R_\fp}^0(\dua{R_\fp}, R_\fp)\cong\Hom_R(\HH_0(\dua{R_\fp}), R_\fp)
\] Since $\HH_0(\dua{R_\fp})\ne 0$ by \ref{size of dc} and $\depth
R_\fp=0$, the module $\Hom_R(\HH_0(\dua{R_\fp}), R_\fp)$ is nonzero,
so we conclude $\Ext^i_R(\dua R,R)_\fp\ne 0$, and hence $\Ext^i_R(\dua
R,R)\ne 0$.
\end{proof}

\begin{corollary}
\label{Cohen-Macaulay}
If $\Ext^i_R(\dua R, R)=0$ for all $i\in[1,\dim R]$, then the ring
$R$ is Cohen-Macaulay.
\end{corollary}

\begin{proof}
Choose $\fp\in\Spec R$ with $\dim (R/\fp)=\dim R$.  Note that
$\fp$ is associated to $R$ and that $\dim(R/\fp)-\depth R=i$, so
$\Ext^i_R(D^R,R)\ne 0$ by the proposition.  In view of our hypothesis,
this implies $i=0$, that is, $R$ is Cohen-Macaulay.
 \end{proof}

A special case of the property \ref{size of dc} reads:

\begin{chunk}
\label{dualizing module}
The ring $R$ is Cohen-Macaulay if and only if $\HH_i(\dua R)=0$ for
all $i\ne 0$.  When this is the case, $\dua R$ can be taken to
be a module, called a {\em dualizing module} of $R$.
\end{chunk}

For the properties of dualizing modules listed below we refer to
\cite[\S 3.3]{BH}.

 \begin{chunk}
 Let $R$ be a Cohen-Macaulay local ring.

\begin{subchunk}
  \label{existence of dualizing}
The ring $R$ has a dualizing module $D$ if and only if there exists a
finite homomorphism of rings $Q\to R$ with a Gorenstein local ring
$Q$; when this is the case, every dualizing module is isomorphic
to $\Ext^{\dim Q-\dim R}_Q(R,Q)$.
 \end{subchunk}

 \begin{subchunk}
 \label{injdim of dualizing}
$\id_RD=\dim R$.
 \end{subchunk}

 \begin{subchunk}
 \label{type}
$\rank_k(D/\fm D)=\rank_k(\Ext^{\dim R}_R(k,R))$.
 \end{subchunk}

 \begin{subchunk}
 \label{ass of dualizing}
$\Ass_R D=\Ass R$.
 \end{subchunk}

 \begin{subchunk}
 \label{localized dualizing}
For each $\fp\in\Spec R$ the module $D_\fp$ is dualizing for $R_\fp$.
 \end{subchunk}

 \begin{subchunk}
 \label{dualfactor}
Each $R$-regular element $f$ is also $D$-regular, and $D/f D$ is dualizing
for $R/(f)$.
 \end{subchunk}

 \begin{subchunk}
 \label{dualizing for gor}
$D\cong R$ if and only if the ring $R$ is Gorenstein.
 \end{subchunk}
 \end{chunk}

\section{Generically Gorenstein rings}
\label{Generically}

A commutative ring $R$ is said to  be {\em generically Gorenstein\/} if the
ring $R_\fp$ is Gorenstein for all $\fp\in \Ass R$.  Our purpose in this
section is to prove

 \begin{theorem}
 \label{CM}
Let $R$ be a generically Gorenstein  local ring.

If $\Ext^i_R(\dua R,R)=0$ for all $i\in[1,\dim R]$, then $R$ is Gorenstein.
 \end{theorem}

 \begin{Remark}
Another proof of Theorem \ref{CM} was independently obtained by Huneke
and Hanes \cite[2.2]{HH} when $R$ is Cohen-Macaulay.  The essential
input comes from \cite[2.1]{HH}, which implies that $\dua
R\otimes_R\dua R$ is a maximal Cohen-Macaulay module.  The same
conclusion can be obtained from Corollary \ref{vs}(3), that has weaker
hypotheses.
 \end{Remark}

We start with a criterion for a module to be free of rank $1$.  With a
view towards later applications, we work in greater generality than
needed here.

 \begin{lemma}
 \label{symmetric}
Let $P\to Q$ be a finite homomorphism of commutative rings, and let $N$
be a finite $Q$-module.  If $P$ is local, $\Ass_P(N\otimes_QN)\cup\Ass_PQ
\subseteq\Ass P$, and $N_\fp\cong Q_\fp$ as $Q_\fp$-modules for
every prime ideal $\fp\in\Ass P$, then $N\cong Q$.
 \end{lemma}

 \begin{proof}
The second symmetric power $\Sym^2_Q(N)$ appears in a  canonical
epimorphism
 \[
\pi^N_Q\col N\otimes_QN\lra \Sym^2_Q(N)
 \]
For each $\fp\in\Ass P$ the homomorphism $(\pi^N_Q)_\fp$ factors as a
composition
 \[
(N\otimes_QN)_\fp\cong N_\fp\otimes_{Q_\fp}N_\fp\xra{\pi^{N_\fp}_{Q_\fp}}
\Sym^2_{Q_\fp}(N_\fp)\cong(\Sym^2_Q(N))_\fp
 \]
with canonical isomorphisms.  As the $Q_\fp$-module
$N_\fp$ is cyclic, the map $\pi^{N_\fp}_{Q_\fp}$ is bijective, so
$(\pi^N_Q)_\fp$ is an isomorphism, hence $(\Ker(\pi^N_Q))_\fp=0$.  The
inclusions
 \[
\Ass_P(\Ker(\pi^N_Q))\subseteq\Ass_P(N\otimes_QN)\subseteq\Ass P
 \]
now imply $\Ker(\pi^N_Q)=0$, so $\pi^N_Q$ is bijective.

Let $\fn$ be a maximal ideal of $Q$ and set $r=\rank_{Q/\fn}(N/\fn
N)$.  Assuming $r=0$, and localizing at $\fn$, we get $N_\fn=0$ by
Nakayama's Lemma.  Since $Q$ is finite as a $P$-module, $\fn\cap P$ is
the maximal ideal of $P$, hence $\fn$ contains $\fp\in\Ass P$.  Thus,
$N_\fp$ is a localization of $N_\fn$, hence $N_\fp=0$, contradicting
our hypothesis.  We conclude that $r\ge1$.  By naturality, the map
$\pi^N_Q$ induces an isomorphism
 \[
(N/\fn N)\otimes_{Q/\fn}(N/\fn N)\cong \Sym^2_{Q/\fn}(N/\fn N)
 \]
of vector spaces over $Q/\fn$.  Comparing ranks, we get $r=1$.

The ring $Q$, being a finite module over a local ring, is semilocal.
Thus, if $\fr$ is its Jacobson radical, then $Q/\fr$ is a finite product
of fields.  By what we have just proved, there is an isomorphism $Q/\fr
Q\cong N/\fr N$ of $Q/\fr Q$-modules.  Nakayama's Lemma implies that the
$Q$-module $N$ is cyclic.  Choose an epimorphism $\varkappa\col Q\to N$.
For each $\fp\in\Ass P$ we have $N_\fp\cong Q_\fp$ by hypothesis, so
$\varkappa_\fp$ is bijective, hence $(\Ker(\varkappa))_\fp=0$.  Since
$\Ass_P(\Ker(\varkappa))$ is contained in $\Ass_PQ$, and by hypothesis
the latter set lies in $\Ass P$, it follows that $\Ker(\varkappa)=0$,
so $\varkappa$ is an isomorphism.
 \end{proof}

 \begin{proof}[Proof of Theorem {\em\ref{CM}}]
By Corollary \ref{Cohen-Macaulay} the ring $R$ is Cohen Macaulay, so in view
of \ref{dualizing module} we may assume that $\dua R$ is an $R$-module, $D$.
Since $\id_RD=\dim R$, cf.\ \ref{injdim of dualizing}, Corollary \ref{vs}
yields an isomorphism $D\otimes_RD\cong\Hom_R(D^*,D)$. This gives the
first equality in the chain
 \[
\Ass_R(D\otimes_RD)=\Ass_R\Hom_R(D^*,D)=\Supp_R(D^*)\cap\Ass_RD \subseteq
\Ass R
 \]
where the second equality comes from a classical expression for the
associator of a module of homomorphisms, and the last inclusion from
property \ref{ass of dualizing}.

As $R_\fp$ is Gorenstein for each $\fp\in\Ass R$, we have
$D_\fp\cong R_\fp$ by \ref{localized dualizing} and \ref{dualizing
for gor}.

Lemma \ref{symmetric} now yields $D\cong R$, so $R$ is Gorenstein by
\ref{dualizing for gor}.
 \end{proof}

\section{Homomorphic images of generically Gorenstein rings}
\label{Residue}

In this section, $Q$ denotes a generically Gorenstein
homomorphic image of a Gorenstein local ring.

Admitting slightly more vanishing, we extend the scope of Theorem \ref{CM}.

 \begin{theorem}
 \label{ci}
Assume $R\cong Q/(\bsg)$, where $\bsg$ is a $Q$-regular set.

If $\Ext^i_R(\dua R,R)=0$ for all $i\in[1,\dim R+1]$, then $R$ is
Gorenstein.
 \end{theorem}

We start with a lemma that grew out of a discussion with Bernd Ulrich.

 \begin{lemma}
  \label{generic}
If $Q$ is Cohen-Macaulay, then $R\cong P/(g)$, where $P$ is a
generically Gorenstein, Cohen-Macaulay homomorphic image of a
Gorenstein local ring and $g$ is $P$-regular.
 \end{lemma}

  \begin{proof}
Being a homomorphic image of a Gorenstein ring, $Q$ contains an ideal
$\fd$ such that for every $\fq\in\Spec Q$, the ring $Q_\fq$ is
Gorenstein if and only if $\fq\supseteq\fd$.  Let $\fn$ denote the
maximal ideal of $Q$ and set $X=\{\fq\in\Spec Q\,|\,\height\fq=1\text{
and }\fq\supseteq\fd\}$.  If $\card(\bsg)=s>1$ then $X$ is finite, so
pick $g_1\in(\bsg)\smallsetminus\big(\fn(\bsg)\cup\bigcup_{\fq\in
X}\fq\big)$ using prime avoidance.  The ideal $\fd'=\fd+(g_1)/(g_1)$
defines the non-Gorenstein locus of $Q'=Q/(g_1)$ and has
$\height\fd'>0$, so $Q'$ is generically Gorenstein.  Extend $g_1$ to a
minimal set $g_1,\dots,g_s$ generating $(\bsg)$.  This set is
$Q$-regular.  Thus, $Q'$ is a Cohen-Macaulay homomorphic image of a
Gorenstein local ring, the images of $g_2,\dots,g_s$ in $Q'$ form a
$Q'$-regular set $\bsg'$, and $Q'/(\bsg')\cong R$.  To obtain the
desired result, iterate the procedure $s-1$ times.
 \end{proof}

Next we record an easy result on change of rings.

 \begin{lemma}
 \label{reduction}
Let $P$ be a local noetherian ring, let $M$, $N$ be finite $P$-modules,
let $g\in P$ be a $P\oplus M\oplus N$-regular element, and let $n$ be
a positive integer.
 \begin{enumerate}[\quad\rm(1)]
\item
If $\Ext_{P/(g)}^n(M/gM,N/gN)=0$, then $\Ext_P^n(M,N)=0$.
\item
If $\Ext^i_P(M,N)=0$ for $i=n,n+1$, then $\Ext^n_{P/(g)}(M/gM,N/gN)=0$.
 \end{enumerate}
 \end{lemma}

 \begin{proof}
The canonical isomorphisms $\Ext_P^{i+1}(M/gM, N)\cong
\Ext_{P/(g)}^{i}(M/gM,N/gN)$, holding for all $i\ge 0$, show that
the exact sequence of $P$-modules
 \[
0\to M \xra{\ g\ }M\to M/gM\to 0
 \]
yields for every $i\ge0$ an exact sequence
 \[
\Ext_P^i(M,N)\xra{\ g\ }\Ext_P^i(M,N)\lra
\Ext_{P/(g)}^{i}(M/gM,N/gN)\lra\Ext_P^{i+1}(M,N)
 \]
This sequence establishes (2).  It also shows that the hypothesis of
(1) implies $g\Ext_P^i(M,N)=0$, so Nakayama's Lemma yields the desired
assertion.
 \end{proof}

 \begin{proof}[Proof of Theorem {\em\ref{ci}}]
By Corollary \ref{Cohen-Macaulay} the ring $R$ is Cohen Macaulay,
so Lemma \ref{generic} yields $R\cong P/(g)$, with $g$ a $P$-regular
element.  From Lemma \ref{reduction}(1) and \ref{dualfactor} we then
get $\Ext^i_P(\dua P,P)=0$ for all $i\in[1,\dim P]$.  Thus, the ring $P$
is Gorenstein by Theorem \ref{CM}, hence so is $R$.
 \end{proof}

\section{Well linked Cohen-Macaulay rings}
\label{Well}

In this section we describe a class of rings to which Theorem \ref{ci}
applies.  They are constructed by means of linkage, so we start by
recalling some terminology.

Let $Q$ be a Gorenstein local ring and $\fa$ an ideal in $Q$ such
that the ring $Q/\fa$ is Cohen-Macaulay.  An ideal $\fb$ of $Q$ is said
to be {\em linked\/} to $\fa$ if $\fb=(\bsg):\fa$ and $\fa=(\bsg):\fb$
for some $Q$-regular sequence $\bsg$ in $\fa\cap\fb$; when this is the
case, we write $\fb\sim\fa$; the ring $Q/\fb$ is Cohen-Macaulay, cf.\
Peskine and Szpiro \cite[1.3]{PS}.

We say that an ideal $\fb$ is {\em generically complete intersection\/},
if for every $\fq\in\Ass_Q(Q/\fb)$ the ideal $\fb_\fq$ in $Q_\fq$ is
generated by a $Q_\fq$-regular sequence.

 \begin{theorem}
 \label{licci}
Let $R$ be a Cohen-Macaulay local ring of the form $Q/\fa$, where $Q$
is a Gorenstein local ring and $\fa$ is an ideal for which there is a
sequence of links $\fa\sim\fb_1\sim\cdots\sim\fb_s\sim\fb$ with $\fb$
a generically complete intersection ideal.

If $\Ext^i_R(\dua R,R)=0$ for all $i\in[1,\dim R+1]$, then $R$ is
Gorenstein.
 \end{theorem}

The proof uses Theorem \ref{ci}, as well as work of Huneke and Ulrich
\cite{HuUl1}, \cite{HuUl2}.

 \begin{proof}
The construction of generic links of $\fb$ in \cite[2.17(a)]{HuUl2}
provides a prime ideal $\fr$ in some polynomial ring $Q[\bsx]$, an
ideal $\fc$ in the local ring $Q[\bsx]_{\fr}$, a regular set $\bsg'$
in $Q'= Q[\bsx]_{\fr}/\fc$, and an isomorphism $R\cong Q'/(\bsg')$.
The ring $Q'$ is Cohen-Macaulay because $R$ is.  By \cite[2.9(b)]{HuUl1}
the ideal $\fc$ is generically complete intersection along with $\fb$.
This implies that $Q'$ is generically Gorenstein.  It is also a residue
ring of the Gorenstein ring $Q[\bsx]_{\fr}$, so Theorem \ref{ci} implies
that the ring $R$ is Gorenstein.
 \end{proof}

The theorem covers the case when $R$ is in the {\em linkage class
of a complete intersection\/}, that is, in some Cohen presentation
$\wh R\cong Q/\fa$ with $Q$ a regular local ring, $\fa$ is linked to
an ideal generated by a $Q$-regular sequence.  Cohen-Macaulay rings
of codimension at most $2$ are of this type, cf.\ \cite[3.3]{PS},
so we obtain:

 \begin{corollary}
 \label{codim 2}
Let $R$ be a Cohen-Macaulay local ring with $\codepth R\le 2$.

If $\Ext^i_R(\dua R,R)=0$ for all $i\in[1,\dim R+1]$, then $R$ is
Gorenstein.
 \qed
 \end{corollary}

\section{Cohen-Macaulay rings with short reductions}
 \label{Loewy}

In this section we consider higher dimensional versions of artinian
rings whose maximal ideal has a low degree of nilpotence.

 \begin{theorem}
\label{loewy}
Let $(R,\fm,k)$ be a local ring containing an $R$-regular sequence
$\bsf$ with the property that  $\fm^3\subseteq(\bsf)$.

If $\Ext^i_R(\dua R,R)=0$ for all $i\in[1,\dim R+1]$, then $R$ is
Gorenstein.
 \end{theorem}

\begin{Remark}
In the special case when $R$ is a $k$-algebra and $\rank_k(R)$ is finite,
this is due to Asashiba \cite{As}; a simplified proof is given by Asashiba
and Hoshino \cite{AsHo}.  A close reading of that proof shows that the
hypothesis that $R$ is a $k$-algebra can be avoided.  We present that
argument, referring to \cite{As} and \cite{AsHo} whenever possible, and
sketching modifications when necessary.  A different proof of the theorem
was obtained by Huneke, \c Sega, and Vraciu \cite{HSV}.  \end{Remark}

 \begin{proof}
By Corollary \ref{Cohen-Macaulay} the ring $R$ is Cohen Macaulay,
so in view of \ref{dualizing module} we may assume that $\dua R$ is an
$R$-module, $D$.  We may further assume $\fm^3=0$ and $\Ext^1_R(D,R)=0$.
Indeed, we get $\Ext^1_{R/(\bsf)}({D/\bsf D},{R/(\bsf)})=0$ by repeated
applications of Lemma \ref{reduction}.  Now note that $D/\bsf D$ is
a dualizing module for $R/(\bsf)$ by \ref{dualfactor}, and that $R$
and ${R/(\bsf)}$ are simultaneously Gorenstein.

For the rest of the proof we fix a free cover
 \begin{equation}
  \label{kernel}
0\lra C \lra F\lra D \lra 0
 \end{equation}

 \begin{step}[{\cite[2.2]{As}}]
The module $C$ is indecomposable.
 \end{step}

 \begin{step}[{\cite[2.3]{As}}]
  \label{socle}
There is an equality $\fm^2=(0:\fm)_R$.
 \end{step}

The proofs of the two steps above do not use the hypothesis that
$R$ is an algebra.

If $M$ is an $R$-module, then we let $\ell(M)$ denote its length.
The polynomial $[M]=\sum_{i\ges0}\rank_k(\fm^i M/\fm^{i+1}M)t^i$ is called
the {\em Hilbert series\/} of $M$.

 \begin{step}[{\cite[2.4 and its proof]{As}}]
  \label{small socle}
There is an inequality $\ell(\fm^2)\le2$.

If equality holds, then $\ell(C/\fm C)=\ell(\fm/\fm^2)$.
 \end{step}

The original proof of Step \ref{small socle} proceeds through computations
with Hilbert series, using Steps \ref{socle} and \ref{small socle}
for rank counts in short exact sequences of $k$-vector spaces; one only
needs to replace the latter by counts of lengths of $R$-modules.

The assertion of the following step is proved without restrictions on $R$.

 \begin{step}[{\cite[2.1]{AsHo}}]
If $\ell(D/\fm D)=2$, then there is an exact sequence
 \begin{equation}
  \label{exact}
0\lra\Tor^R_1(D,D)\lra C\otimes_RD\lra\Hom_R(D,D)
 \end{equation}
 \end{step}

If $\ell(M)$ is finite, then for every $x\in R$ a length count
in the exact sequence
 \[
0\lra(0:x)_M\lra M\xra{x}M\lra M/xM\lra 0
 \]
yields $\ell((0:x)_M)=\ell(M/xM)$.  Such equalities are used in the
proof below.

 \begin{step}[{\cite[Proofs of 2.2--2.4 and 3.4]{AsHo}}]
The ring $R$ is Gorenstein.
 \end{step}

If $\fm^2=0$, then $\fm C=0$, so the isomorphism
$\Ext^1_R(D,R)\cong\Ext^2_R(C,R)$ induced by \eqref{kernel} yields
$\Ext^2_R(k,R)=0$, and hence $R$ is Gorenstein.

If $\ell(\fm^2)=1$, then $(0:\fm)\cong k$ by Step \ref{socle}, so
$R$ is Gorenstein.

By Step \ref{small socle}, for the rest of the proof we may assume
$\ell(\fm^2)=2$ and $\ell(C/\fm C)=\ell(\fm/\fm^2)$.  From \ref{type}
and Step \ref{socle} we get $\ell(D/\fm D)=2$.  Since $D$ is a faithfully
injective $R$-module and $\Hom_R(D,D)\cong R$, we have
 \[
\Hom_R(\Tor^R_1(D,D),D)\cong\Ext^1_R(D,\Hom_R(D,D))\cong\Ext^1_R(D,R)
 \]
and hence $\Tor^R_1(D,D)=0$.   Thus, from \eqref{exact} we get an
injection $C\otimes_RD\hookrightarrow R$.  Since $\fm^2\ne0$, we
may choose an element $x\in\fm\smallsetminus(0:\fm)_R$.  We now have
 \begin{align*}
1+\ell(\fm/\fm^2)
&=\ell(R)-2\\
&\ge\ell(R)-\ell(Rx)
=\ell\big((0:x)_R\big)\\
&\ge\ell\big((0:x)_{C\otimes_RD}\big)
=\ell\big((C\otimes_RD)/x(C\otimes_RD)\big)\\
&\ge\ell\big((C\otimes_RD)/\fm(C\otimes_RD)\big)
=\ell(C/\fm C)\ell(D/\fm D)\\
&=2\ell(\fm/\fm^2)
 \end{align*}
As a consequence, we get $\ell(\fm/\fm^2)\le1$,
hence $\ell(\fm^2)\le1$, a contradiction.
 \end{proof}

\section{Artinian rings with Gorenstein deformations}
 \label{Artinian}

Given a homomorphism of commutative rings $P\to Q$, let $\du PQ$ denote
the $Q$-module $\Hom_P(Q,P)$ with the canonical action. The $P$-algebra
$Q$ is said to be {\em Frobenius\/} if $Q$ is a finite free $P$-module
and there is a $Q$-linear isomorphism $\du PQ\cong Q$.  For every prime
ideal $\fq$ of $P$ we set $k(\fq)=P_{\fq}/\fq P_{\fq}$.

 \begin{theorem}
 \label{deformation}
Let $(P,\fp,k)$ be a local ring and let $P\to Q$ be a homomorphism of
commutative rings, such that $Q$ is a finite free $P$-module, and there
exists a prime ideal $\fq$ of $P$ for which the ring $k(\fq)\otimes_{P}Q$
is Gorenstein.

If $R=Q/\fp Q$ satisfies $\Ext^{1}_{R}(\du kR,R) = 0$, then the
ring $R$ is Gorenstein and the $P$-algebra $Q$ is Frobenius.
 \end{theorem}

 \begin{remark}
  \label{products}
The ring $R$ in the theorem is a finite $k$-algebra.  For such rings,
$X=\Spec R$ is finite and $R\cong\prod_{\fm\in X}R_\fm$ as $k$-algebras.
This yields the first isomorphism of $R$-modules in the next formula,
the second one comes from \ref{existence of dualizing}:
 \[
\du kR\cong\coprod_{\fm\in X}\du k{R_\fm}\cong\coprod_{\fm\in
X}\dua{R_\fm}
 \]
  \end{remark}

For use in the proof we recall a couple of known facts, cf.\ e.g.\
\cite[\S 14]{SS}.

 \begin{chunk}
 \label{decomposition}
Let $P\to P'$ be a homomorphism of commutative rings, set
$Q'=P'\otimes_PQ$, and let $P'\to Q'$ be the induced homomorphism.
In the composition
 \[
Q'\otimes_Q\du PQ\lra P'\otimes_P\du PQ\lra\du {P'}{Q'}
 \]
of canonical homomorphisms of $Q'$-modules the first arrow is always
bijective.  The second one is an isomorphism when $Q$ is finite free
over $P$.
 \end{chunk}

 \begin{lemma}
 \label{Frobenius}
If $(P,\fp,k)$ is a local ring and $P\to Q$ a homomorphism
of commutative rings, then the following are equivalent.
 \begin{enumerate}[\quad\rm(i)]
  \item
The $P$-algebra $Q$ is Frobenius.
 \item
The $P'$-algebra $Q'=P'\otimes_PQ$ is Frobenius for every local ring $P'$
and every homomomorphism of rings $P\to P'$.
 \item
The $P$-module $Q$ is finite free and the ring $k(\fq)\otimes_PQ$ is
Gorenstein for every $\fq\in\Spec P$.
 \item
The $P$-module $Q$ is finite free and the ring $R=Q/\fp Q$ is Gorenstein.
 \end{enumerate}
  \end{lemma}

 \begin{proof}
Consider first the special case $P=k$.  The isomorphisms of Remark
\ref{products} then show that the ring $R$ is Gorenstein if and only
if $R_\fm$ is Gorenstein for every $\fm$, and the $k$-algebra $R$ is
Frobenius if and only if $\du k{R_\fm}\cong R_\fm$ for every $\fm$.
The conditions at each $\fm$ are equivalent by \ref{dualizing for
gor}, so (i) $\iff$ (iv) whenever $P$ is a field.

Under the hypothesis of (i) the second map in \ref{decomposition}
is bijective, so (ii) holds.  When (ii) holds the algebra $k(\fq)\otimes_PQ$
over the field $k(\fq)$ is Frobenius.  By the
special case this implies that the ring $k(\fq)\otimes_PQ$ is Gorenstein,
which is the assertion of (iii). It is clear that (iii) implies (iv).
If (iv) holds, then from Remark \ref{products} and \ref{dualizing for
gor} we obtain an isomorphism of $R$-modules $k\otimes_Q\du PQ\cong R$.
Nakayama's Lemma then yields an epimorphism of $Q$-modules $Q\to\du PQ$.
Over $P$ both modules are finite free and have the same rank, so this
map is bijective.
 \end{proof}

 \begin{proof}[Proof of Theorem \rm\ref{deformation}]
By Lemma \ref{Frobenius}, the $P$-algebra $Q$ is Frobenius if and only
if the ring $R$ is Gorenstein.  The induced map $P/\fq\to Q/\fq Q$ turns
$Q/\fq Q$ into a finite free $P/\fq$-module; due to the isomorphism
of rings $(Q/\fq Q)/\fp(Q/\fq Q)\cong R$, the same lemma implies that $R$
is Gorenstein if and only the $P/\fq$-algebra $Q/\fq Q$ is Frobenius.

Changing notation if necessary, for the rest of the proof we may
assume that $P$ is a domain with field of fractions $K$, and the ring
$K\otimes_PQ$ is Gorenstein.  Setting $D=\du PQ$, we now have isomorphisms
of $K\otimes_PQ$-modules
\[
K\otimes_PD\cong\du K{K\otimes_PQ}\cong K\otimes_PQ
\]
given by \ref{decomposition} and Lemma \ref{Frobenius}.  In view of
Lemma \ref{symmetric}, the desired isomorphism $D\cong Q$ will
follow once we prove that the $P$-module $D\otimes_QD$ is torsion-free.

Let $F\xra{\simeq}D$ be a resolution by free $Q$-modules, and let
$k\xra{\simeq}J$ be a resolution by injective $P$-modules.  There are
then chains of relations
  \begin{gather}
 \begin{aligned}
  \label{first chain}
R\otimes_{Q}F
&\cong(k\otimes_PQ)\otimes_QF
\cong k\otimes_PF
\simeq k\otimes_PD\\
&\cong\du kR
 \end{aligned}
\\
 \begin{aligned}
  \label{second chain}
\Hom_{P}(D,J)
&\simeq\Hom_{P}(D,k)
\cong\Hom_{k}(k\otimes_{P}D,k)
\cong\Hom_k(\du kR,k)\\
&\cong R
 \end{aligned}
  \end{gather}
where the quasi-isomorphisms result from the freeness of $D$ over
$P$, cf.\ \eqref{semiflat} and \eqref{semiprojective}, the isomorphism
$k\otimes_PD\cong\du kR$ comes from \ref{decomposition},
and all other maps are canonical.  Using \eqref{second chain} and
\eqref{semiprojective} we get quasi-isomorphisms
 \begin{align*}
\Hom_{R}(R\otimes_QF,R)
&\cong\Hom_{Q}(F,R)
\simeq\Hom_{Q}(F,\Hom_{P}(D,J))\\
&\cong\Hom_{P}(F\otimes_{Q}D,J)
 \end{align*}
As \eqref{first chain} shows that $R\otimes_{Q}F$ is a free resolution of
the $R$-module $\du kR$, Proposition \ref{spectral sequence} applied to
$G=F\otimes_{Q}D$ and $J$ yields a strongly convergent spectral sequence
 \[
\EE2pq=\Ext_{P}^{-p}(\Tor_{-q}^{Q}(D,D),k)\implies\Ext^{-p-q}_R(\du kR,R)
 \]
with differentials $\DD rpq\col\EE rpq\to \EE r{p+r}{q+r-1}$.  Since $\EE
2pq=0$ if $p>0$ or $q>0$, the spectral sequence defines an exact sequence
 \[
0\lra\Ext^{1}_{P}(D\otimes_{Q}D,k)\lra\Ext^{1}_R(\du kR,R)
 \]
which implies that $\Ext^{1}_{P}(D\otimes_{Q}D,k)$ vanishes.  Because the
$P$-module $D\otimes_{Q}D$ is finite, we conclude that it is actually
free, and so {\em a fortiori\/} torsion-free.
 \end{proof}

\section{Rings of small codepth}
\label{Codepth}

In this section we prove the theorem below, after substantial preparation.

 \begin{theorem}
  \label{SC}
Let $R$ be a local ring with $\codepth R\le 3$.

If $\Ext^i_R(\dua R,R)=0$ for all $i\gg 0$, then $R$ is Gorenstein.
 \end{theorem}

Let $M$, $N$ be complexes with finite homology.

\setcounter{theorem}{2}
 \begin{subchunk}
  \label{finiteext}
The $R$-modules $\Tor^R_i(M,N)$ and $\Ext_R^i(M,N)$ are finite
for each $n\in\BZ$, and vanish for all $n\ll0$.
 \end{subchunk}

The preceding result is established by using standard arguments with
spectral sequences.  It ensures that the right hand sides of the equalities
 \begin{align*}
\Ba_R^M(t)&=\sum_{n\in\BZ}\rank_k\big(\Ext_R^i(k,M)\big)t^i
 \\
\Po^R_M(t)&=\sum_{n\in\BZ}\rank_k\big(\Tor^R_i(M,k)\big)t^i
 \end{align*}
are formal Laurent series.  They are called, respectively, the {\em Bass
series\/} and the {\em Poincar\'e series\/} of $M$, and are invariant
under quasi-isomorphisms by \ref{extcong}, \ref{torcong}.

 \begin{subchunk}
 \label{truncation}
If $F\to M$ is a semiprojective resolution and $m=\sup\HH(M)$,
then the module $M'=\Coker(\dd^F_{m+1})$ satisfies
$\Po^R_M(t)-t^{m}\Po^R_{M'}(t)\in\BZ[t,t^{-1}]$.

Indeed, $\Tor^R_i(M,k)\cong\Tor^R_{i-m}(M',k)$ holds for each $i>m$ by
Lemma \ref{shifts}.
 \end{subchunk}

The next equality is due to Foxby, cf.\ \cite[4.1(a)]{Foxby-iso} and
\cite[15.18]{Foxby-notes}.

 \begin{subchunk}
  \label{formula}
$\Ba^{\rh_R(M,N)}_R(t)=\Po_M^R(t)\Ba^N_R(t)$
 \end{subchunk}

Once again, we focus upon the special case of dualizing complexes.

 \begin{chunk}
Set $d=\depth R$ and $D=\dua R$.

 \begin{subchunk}
  \label{dualba}
$\Ba^{D}_R(t)=t^{d}$ holds by \ref{id of dc}.
 \end{subchunk}

We recall another fundamental property of dualizing complexes.

 \begin{subchunk}
  \label{duality}
Set $M^{\dagger}=\rh_R(M,D)$.  The canonical map $M\to M^{\dagger\dagger}$
is then a quasi-isomorphism, cf.\ \cite[15.14]{Foxby-notes}.  For $M=R$
this yields $R\simeq\rh_R(D,D)$.
  \end{subchunk}

\setcounter{equation}{2}
 From the formulas and quasi-isomorphisms above one now gets:
 \begin{gather}
 \label{dualpo}
\Ba_R^R(t)=\Ba_R^{\rh_R(D,D)}(t)=\Po_D^R(t)\Ba^D_R(t)=
\Po_D^R(t)t^{d}\\
 \label{bapo}
\Ba_R^{M^\dagger}(t)=\Po_M^R(t)\Ba^D_R(t)=\Po_M^R(t)t^{d}\\
 \label{poba}
\Ba_R^M(t)=\Ba_R^{M^{\dagger\dagger}}(t)=
\Po_{M^\dagger}^R(t)\Ba^D_R(t)=\Po_{M^\dagger}^R(t)t^{d}
 \end{gather}
  \end{chunk}

For the proof the theorem we need specific information on Poincar\'e
series and Bass series of complexes.  It suffices to have it for
Poincar\'e series of modules:

 \begin{proposition}
  \label{equivalences}
If $R$ is a local ring that has a dualizing complex and $c(t)$ is a
polynomial in $\BZ[t]$, then the following statements are equivalent:
  \begin{enumerate}[\quad\rm(1)]
 \item
$c(t)\Po^R_M(t)\in\BZ[t]$ for all finite $R$-modules $M$.
 \item
$c(t)\Po^R_M(t)\in\BZ[t,t^{-1}]$ for all complexes $M$ with finite
homology.
 \item
$c(t)\Ba^M_R(t)\in\BZ[t,t^{-1}]$ for all complexes $M$ with finite
homology.
 \end{enumerate}
  \end{proposition}

 \begin{proof}
Formulas \eqref{bapo} and \eqref{poba} show that (2) and (3) are equivalent.
Obviously, (2) implies (1), and the converse follows from \ref{truncation}.
 \end{proof}

We need several results on Poincar\'e series of modules
over rings of small codepth.

 \begin{chunk}
Let $(R,\fm,k)$ be a ring, such that $\codim R\le 3$.

The next result is proved by Avramov, Kustin and Miller \cite[6.1 and
5.18]{AKM}.

 \begin{subchunk}
  \label{rational}
There exists a polynomial $d(t)\in\BZ[t]$, such that
$(1+t)d(t)\Po_M^R(t)\in\BZ[t]$ for every finite $R$-module $M$ and
$(1+t)d(t)\Po_k^R(t)=(1+t)^{\edim R}$.
 \end{subchunk}

The table below gives a classification of the rings of codepth at most
$3$ obtained in \cite[2.1]{AKM}; the values of $d(t)$ are computed by
Avramov \cite[3.5]{Av1}.

 \begin{subchunk}
 \label{table}
If $R$ is not complete intersection, then there exist integers $l$, $m$,
$p$, $q$, $r$ that determine the polynomial $d(t)$ from \ref{rational}
through the following table:

\bigskip
  \begin{center}
\renewcommand{\arraystretch}{1.5}
   \begin{tabular}{|c|c|c|c|}
   \hline
    type       &$\codepth R$  &$d(t)$                      & restrictions
\\ \hline \hline
    $\GO$      &2           &$1-t-lt^2$                  & $l\ge 1$
\\ \hline
    $\TE$      &3           & $1-t-lt^2-(m-l-3)t^3-t^5$  & $m>l+1\ge 3$
\\ \hline
    $\B$       &3           & $1-t-lt^2-(m-l-1)t^3+t^4$  & $m>l+1\ge 3$
\\ \hline
    $\G(r)$    &3           & $1-t-lt^2-(m-l)t^3+t^4$    & \begin{tabular}{c}
                                                            $m>l+1\ge 3$
                                                         \\
                                                           $l+1\ge r\ge 2$
                                                           \end{tabular}
\\ \hline
    $\Ho(p,q)$ &3           &$1-t-lt^2-(m-l-p)t^3+qt^4$  & \begin{tabular}{c}
                                                            $m>l+1\ge 3$
                                                         \\
                                                           $l\ge p\ge 0$
                                                         \\
                                                           $m-l\ge q\ge 0$
                                                           \end{tabular}
\\ \hline
   \end{tabular}
 \end{center}
 \end{subchunk}

As a consequence of \cite[3.6]{Av1}, we have:

 \begin{subchunk}
  \label{pole}
If $d(t)$ has $1$ as a root, then $R$ is of type $\Ho(l,l-1)$ for some
$l\ge2$ and $d(t)=(1+t)(1-t)(1-t-(l-1)t^2)$.
 \end{subchunk}

Sun \cite[Proof of 1.2, p.\ 61]{Su} provides additional information on the
roots of $d(t)$.

 \begin{subchunk}
  \label{Sun}
If $d(t)$ has a real root $r$ with $0<r<1$, then this root is simple.
 \end{subchunk}
 \end{chunk}

In the proof of the theorem we use a technique developed by \c Sega \cite{Se}.
In order to apply it, we need an extension of Sun's result.

 \begin{lemma}
  \label{multiplicity}
Let $R$ be a local ring with $\codim R\le 3$, which is not complete
intersection, and let $d(t)$ be the polynomial described in
{\rm\ref{rational}}.

If $p(t)\in\BZ[t]$ is an irreducible polynomial with constant term $1$ and
at least one negative coefficient, then $p(t)^2$ does not divide $d(t)$.
 \end{lemma}

\begin{proof}
We assume that our assertion fails, note that it implies $\deg
p(t)\le2$, and obtain a contradiction for each type of rings described
in \ref{table}.

$\GO$.  We must have $1-t-lt^2=p(t)^2$, hence $l=-1/4$.  This is absurd.

$\TE$.  Only $1+t$ can be a linear factor of $d(t)$, so $p(t)=1+at\pm
t^2$ and $d(t)=(1+t)p(t)^2=1+\cdots+t^5$.  However, by \ref{table}
the coefficient of $t^5$ is $-1$.

$\B$ or $\G(r)$.  Only $1+t$ can be a linear factor of $d(t)$.  Thus,
$p(t)=1+at\pm t^2$ and $d(t)=p(t)^2=1+2at+\cdots$, yielding $a=-1/2$,
a contradiction.

$\Ho(p,q)$.  If $p(t)=1+at+bt^2$, then $d(t)=p(t)^2$, so we get a
contradiction as above.  If $p(t)=1-at$ with $a>1$, then $d(t)$ has
a double real root $r=1/a<1$, contradicting \ref{Sun}.  Finally, if
$p(t)=1-t$, then $d(t)=(1+t)(1-t)(1-t-(l-1)t^2)$ by \ref{pole}, so $1-t$
divides $1-t-(l-1)t^2$, hence $l=1$.  This is impossible.
 \end{proof}

\setcounter{theorem}{6}

 \begin{proof}[Proof of Theorem {\rm \ref{SC}}]
Using {\rm\ref{rational}} and \ref{equivalences}, choose
$d(t),r(t)\in\BZ[t]$ such that
 \begin{equation}
  \label{bass}
\Ba_R^R(t)=\frac{r(t)}{(1+t)d(t)}
 \end{equation}
Set $D=\dua R$ and $d=\depth R$.  Our hypothesis means that
$\HH\rh_R(D,R)$ is bounded.  It is degreewise finite by \ref{finiteext},
so using \ref{formula} and \eqref{dualpo} we get
 \[
\Ba_R^{\rh_R(D,R)}(t)=\Po^R_D(t)\Ba_R^R(t)=\frac{\Ba_R^R(t)^2}{t^d}=
\frac{r(t)^2}{t^d(1+t)^2d(t)^2}
 \]
Referring again to \ref{rational} and \ref{equivalences}, we now obtain
 \begin{equation}
  \label{divides}
\frac{r(t)^2}{d(t)}= {t^d}(1+t)\bigg((1+t)d(t)\Ba_R^{\rh_R(D,R)}(t)\bigg)
\in\BZ[t,t^{-1}]
 \end{equation}
Let $s(t)$ denote the product of all the irreducible in $\BZ[t]$ factors
of $d(t)$ with constant term $1$ and at least one negative coefficient,
and set $q(t)=d(t)/s(t)$.

If $R$ is complete intersection, then it is Gorenstein, so there is
nothing to prove.  Else, \eqref{divides} and \ref{multiplicity} imply
that $s(t)$ divides $r(t)$ in $\BZ[t,t^{-1}]$, so \eqref{bass} yields
 \[
(1+t)q(t)\Ba_R^R(t)=\frac{r(t)}{s(t)}\in\BZ[t,t^{-1}]
 \]
Since $q(t)$ and $\Ba_R^R(t)$ have non-negative coefficients, the last
formula implies $\Ba_R^R(t)\in\BZ[t,t^{-1}]$, hence $\id_RR<\infty$.
Thus, $R$ is Gorenstein, as desired.
 \end{proof}

\section{Golod rings}
\label{Golod}

Let $(R,\fm,k)$ be a local ring.  Serre proved a coefficientwise
inequality
 \[
\Po_{k}^R(t)\preccurlyeq\frac{(1+t)^{\edim R}}{1-\sum_{j=1}^{\infty}\rank
\HH_{j}(K^R)t^{j+1}}
 \]
of formal power series, where $K^R$ denotes the Koszul complex on
a minimal set of generators of $\fm$.  If equality holds, then $R$
is said to be a {\em Golod ring\/}.

The class of Golod rings is essentially disjoint from that of Gorenstein
rings: if $R$ belongs to their intersection, then $R$ is a {\em
hypersurface ring\/}, in the sense that its $\fm$-adic completion $\wh R$
is the homomorphic image of a regular local ring by a principal ideal.

 \begin{theorem}
 \label{golod}
Let $R$ be a Golod local ring.

If $\Ext^i_R(\dua R,R)=0$ for all $i\gg 0$, then $R$ is a hypersurface
ring.
 \end{theorem}

For the proof we recall a theorem of Jorgensen \cite[3.1]{Jo}:

 \begin{chunk}
 \label{Jorgensen}
Let $R$ be a Golod ring and $L$, $M$ modules with finite
homology.

If $\Tor^R_i(L,M)=0$ for all $i\gg 0$, then $L$ or $M$ has finite
projective dimension.
 \end{chunk}

We need an extension of this result to complexes.

 \begin{proposition}
 \label{vanishing Tors}
Let $R$ be a Golod ring and $L$, $M$ complexes with finite homology.

If $\Tor^R_i(L,M)=0$ for all $i\gg 0$, then $\Po^R_{L}(t)$ or
$\Po^R_{M}(t)$ is a Laurent polynomial.
 \end{proposition}

 \begin{proof}
Set $l=\sup\HH(L)$ and $m=\sup\HH(M)$, choose semiprojective resolutions $E\to
L$ and $F\to M$, then set $L'=\Coker(\dd^E_l)$ and $M'=\Coker(\dd^F_m)$.
In view of Lemma \ref{shifts}, our hypothesis yields $\Tor_i^R(L',M')=0$
for all $i\gg 0$.  By \ref{Jorgensen}, this means that $\Po^R_{L'}(t)$
or $\Po^R_{M'}(t)$ is a polynomial, so \ref{truncation} gives the
desired assertion.
 \end{proof}

 \begin{proof}[Proof of Theorem {\rm\ref{golod}}]
Set $D=\dua R$ and choose a  semiprojective resolution $F\to D$ with $\inf
F=0$. For all $i\gg 0$ we have $\HH_i\Hom_R(F,R)=0$ by hypothesis,
hence $\HH_{i}(F\otimes_RD)=0$ by Proposition \ref{ss}(2), and thus
$\Po^D_R(t)\in\BZ[t,t^{-1}]$ by Proposition \ref{vanishing Tors}.
Now \eqref{dualpo} implies $\Ba_R^R(t)\in\BZ[t]$, hence
$\id_RR<\infty$, that is, $R$ is Gorenstein.  As noted before the
statement of the theorem, $R$ is then a hypersurface ring.
 \end{proof}

\section{Relations with Tachikawa's conjecture}
\label{Relations}

In this section $A$ denotes a not necessarily commutative algebra of
finite rank over a field $k$, and $\du kA= \Hom_k(A,k)$ has the
canonical action of $A$ on the left.

In \cite{Nak} Nakayama proposed a now famous conjecture.  The original
statement was in terms of the injective resolution of $A$ as a bimodule
over itself, but M\"uller \cite{Mu} proved that it is equivalent to the
one-sided version below:

 \begin{Conjecture}[NC]
If the left $A$-module $A$ has an injective resolution in which each
term is also projective, then $A$ is selfinjective.
 \end{Conjecture}

Considerable efforts notwithstanding, this still is one of the main
open problems in the representation theory of Artin algebras.  It has
generated a significant body of research and has led to the study of
several variants and to the formulation of equivalent forms, see Yamagata
\cite{Ya} for a general presentation of the subject.  In particular,
Tachikawa \cite[Ch.\ 8]{Ta} proved that the validity of Nakayama's
conjecture is equivalent to that of both statements below for all
$k$-algebras of finite rank.

 \begin{Conjecture}[TC1]
If $\Ext^i_A(\du kA,A)=0$ for all $i>0$, then $A$ is selfinjective.
 \end{Conjecture}

 \begin{Conjecture}[TC2]
If the ring $A$ is selfinjective and $M$ is a finite $A$-module
satisfying $\Ext^i_A(M,M)=0$ for all $i>0$, then $M$ is projective.
 \end{Conjecture}

Now we turn to the case when $A$ is commutative.  In view of Remark
\ref{products}, for the conjectures above one may assume that $A$
is local.  The validity of Nakayama's Conjecture is then well-known and
easy to see.  However, both of Tachikawa's Conjectures are open for
commutative algebras.  This apparent paradox is explained by the fact
that the proof of the equivalence of Conjecture (NC) with Conjectures
(TC1) and (TC2) involves intermediate non-commutative algebras.

Here are some instances in which Conjecture (TC1) is known to hold.

 \begin{chunk}
Let $(A,\fr,k)$ be a commutative local algebra of finite rank over $k$.

Each set of hypotheses below implies that the algebra $A$ is
selfinjective:
 \begin{subchunk}
  \label{hoshino}
$\edim A\le 2$ and $\Ext^i_A(\du kA,A)=0$ for $i=1,2$ (Hoshino \cite{Ho}).
 \end{subchunk} \begin{subchunk}
  \label{zeng}
$\rank_k(A/xA)\le 2$ for some  $x\in\fr$ and $\Ext^1_A(\du kA,A)=0$
(Zeng \cite{Ze}).
 \end{subchunk} \begin{subchunk}
  \label{asashiba-hoshino}
$\rank_k(A/xA)\le 3$ for some  $x\in\fr$ and $\Ext^i_A(\du kA,A)=0$
for $i=1,2$ (Asashiba and Hoshino \cite{AsHo}).
 \end{subchunk} \begin{subchunk}
  \label{asashiba}
$\fr^3=0$ and $\Ext^1_A(\du kA,A)=0$ (Asashiba \cite{As}).
 \end{subchunk}
  \end{chunk}

The last result was extended in Theorem \ref{loewy}.  Other results
of this paper establish  Conjecture (TC1) in new cases, some of them
collected in the next theorem.  We remark that Part (1) below improves
\ref{hoshino} and implies \ref{zeng}.

 \begin{theorem}
 \label{artinian}
A commutative local algebra $A$ of finite rank over a field $k$
is selfinjective whenever one of the following
assumptions is satisfied.
 \begin{enumerate}[\quad\rm(1)]
  \item
$\edim A\le 2$ and $\Ext^1_A(\du kA,A)=0$.
  \item
$\edim A\le 3$ and $\Ext^i_A(\du kA,A)=0$ for all $i\gg 0$.
  \item
$A$ is Golod and  $\Ext^i_A(\du kA,A)=0$ for all $i\gg 0$.
  \item
$A\cong Q/(\bsf)$, where $Q$ is a reduced complete local ring,
$\bsf$ is a $Q$-regular sequence, and $\Ext^1_A(\du kA,A)=0$.
 \end{enumerate}
   \end{theorem}

 \begin{proof}
The first three parts are specializations of Corollary \ref{codim 2},
Theorem \ref{SC}, and Theorem \ref{golod}, respectively.  For (4),
use Cohen's Structure Theorem to present $Q$ as a residue ring of
a regular local ring modulo a radical ideal.  Such an ideal
is generically complete intersection, so Theorem \ref{ci} applies.
 \end{proof}

\begin{appendix}

\section{Constructions with complexes}
\label{Constructions}

In this appendix we describe terminology and notation that are used
throughout the paper.  Proofs of the basic results listed below can be
found, for instance, in \cite{Ha}, \cite{Roberts}, or \cite{Foxby-notes}.
We let $R$ denote an associative ring over which modules have left
actions, and $M$ denote a complex of $R$-modules of the form
 \[
\cdots\lra M_{i+1}\xra{\dd^M_{i+1}}M_i \xra{\ \dd^M_i\ }M_{i-1}\lra\cdots
 \]
We write $b\in M$ to indicate that $b$ belongs to $M_i$ for some
$i\in\BZ$, and $|b|=n$ to state that $b$ is in $M_n$.  The $n$'th {\em
shift\/} $\shift^nM$ of $M$ is the complex with
 \[
(\shift^nM)_i=M_{i-n} \qquad\text{and}\qquad
\dd_i^{\shift^nM}=(-1)^n\dd^M_{i-n}\,.
 \]

If $N$ is a complex, then a {\em morphism\/} $\alpha\col M\to N$ is
a sequence of $R$-linear maps $\alpha_i\col M_i\to N_i$,
such that $\dd^N_i\alpha_i=\alpha_{i-1}\dd^M_i$ holds for all $i$.
It is a {\em quasi-isomorphism\/} if $\HH_i(\alpha)$ is bijective for
each $i$; we use the symbol $\simeq$ to identify
quasi-isomorphisms.

\setcounter{theorem}{1}

Set $\inf M=\inf\{i\,|\, M_i\ne 0\}$ and $\sup M=\sup\{i\,|\, M_i\ne 0\}$.
The complex $M$ is {\em bounded above\/} (respectively, {\em below\/})
if $\sup M$ (respectively, $\inf M$) is finite.  It is {\em bounded\/}
when both numbers are finite.  It is {\em degreewise finite} if each
$R$-module $M_i$ is finite.  If $\HH(M)$ is bounded and degreewise finite,
then $M$ has {\em finite homology\/}.

 \begin{subchunk}
  \label{stupid}
We let $M_{<n}$ denote the subcomplex of $M$ having $(M_{<n})_i=0$
for $i\ge n$ and $(M_{<n})_i=M_i$ for $i<n$, and we set $M_{\ges
n}=M/M_{<n}$.
 \end{subchunk}

 \begin{subchunk}
  \label{smart}
We let $\tau_{\les n}(M)$ denote the residue complex
of $M$ with $\tau_{\les n}(M)_i=0$ for $i>n$, $\tau_{\les
n}(M)_n=M_n/\Imag(\dd^M_{n+1})$, and $\tau_{\les n}(M)_i=M_i$ for $i<n$;
the maps $\HH_i(M)\to\HH_i(\tau_{\les n}(M))$ induced by the natural
surjection are bijective for all $i\le n$.
 \end{subchunk}

 \begin{chunk}
Let $L$ be a complex of right $R$-modules and $M$, $N$ be
complexes of $R$-modules.

A complex $\Hom_R(M,N)$ is defined by the formulas
 \[
\Hom_R(M,N)_n=\prod_{j-i=n}\Hom_R(M_i,N_j)\qquad
\dd(\beta)=\dd^N\circ \beta-(-1)^{|\beta|}\beta\circ\dd^M
 \]
and a complex $L\otimes_RM$ defined by the formulas
 \[
(L\otimes_RM)_n=\coprod_{h+i=n}L_h\otimes_RM_i\qquad
\dd(a\otimes b)=\dd^L(a)\otimes b+(-1)^{|a|}a\otimes\dd^M(b)
 \]

If $\lambda\col L\to L'$ is a quasi-isomorphism of complexes of right
$R$-modules, $\mu\col M\to M'$ is a quasi-iso\-morphism of complexes
of $R$-modules, $J$ is a bounded above complex of injective $R$-modules,
and $F$ is a bounded below complex of projective $R$-modules, then the maps
below are quasi-isomorphisms.
 \begin{align}
  \label{semiprojective}
\Hom_R(F,\mu)&\col\Hom_R(F,M)\lra\Hom_R(F,M')\\
  \label{semiinjective}
\Hom_R(\mu, J)&\col\Hom_R(M',J)\lra\Hom_R(M,J)\\
  \label{semiflat}
\lambda\otimes_RF&\col L\otimes_RF\lra L'\otimes_RF
 \end{align}
If, furthermore, $\phi\col F\to F'$ is a quasi-isomorphism of bounded
below complexes of projective $R$-modules, then the maps below are
quasi-isomorphisms.
 \begin{align}
  \label{semiprojective quism}
\Hom_R(\phi,M)&\col \Hom_R(F',M)\lra\Hom_R(F,M)
\\
  \label{semiflat quism}
L\otimes_R\phi&\col L\otimes_RF\lra L\otimes_RF'
 \end{align}
 \end{chunk}

We use derived functors of complexes, for which we recall some
basic notions.

 \begin{chunk}
If $\HH(M)$ is bounded below, then there exists a quasi-isomorphism
$F\to M$ with $F$ a bounded below complex of projective modules,
called a {\em semiprojective resolution\/} of $M$; it can be chosen with
$\inf F=\inf\HH(M)$.  Any two semiprojective resolutions of $M$ are homotopy
equivalent, so the complexes $\Hom_R(F,N)$ and $L\otimes_RF$ are defined
uniquely up to homotopy.  The symbols $\rh_R(M,N)$ and $L\lf_RM$ denote
any complex in the corresponding homotopy class.  One sets
 \begin{align}
  \label{ext}
\Ext^i_R(M,N)&=\HH_{-i}\rh_R(M,N)\\
  \label{tor}
\Tor^R_i(L,N)&=\HH_i(L\lf_RN)
 \end{align}

If $L\to L'$, $M\to M'$ and $N\to N'$ are quasi-isomorphisms, then so are
 \begin{gather*}
\rh_R(M,N)\lra\rh_R(M,N')\longleftarrow\rh_R(M',N')\\ L\lf_RM\lra
L'\lf_RM\lra L'\lf_RM'
 \end{gather*}
by \eqref{semiprojective} and \eqref{semiprojective quism}, respectively
by \eqref{semiprojective} and \eqref{semiprojective quism}.  Thus,
for all $i\in\BZ$
 \begin{align} \label{extcong}
\Ext_R^i(M,N)&\cong\Ext_R^i(M',N')\\
 \label{torcong}
\Tor^R_i(L,M)&\cong\Tor^R_i(L',M')
 \end{align}

\setcounter{subchunk}{4}
 \begin{subchunk}
  \label{modules}
We identify $R$-modules with complexes concentrated in degree zero.
For modules the constructions above yield the classical derived functors.
 \end{subchunk}
  \end{chunk}

\section{Lemmas about complexes}
\label{Lemmas}

The hypotheses here are the same as in Appendix \ref{Constructions}.
In one form or another, the results collected below are probably known
to (some) experts, but may not be easy to find in the literature.
Furthermore, slight weakenings of the hypotheses may lead to false
statements, so we are providing complete arguments.

We start with the computation of a classical spectral sequence.

 \begin{proposition}
  \label{spectral sequence}
If $G$ is a complex of $R$-modules and $J$ is a complex of injective
$R$-modules with $\sup J=0$, then there exists a spectral
sequence having
 \[
\EE 2pq=\HH_p\Hom_R(\HH_{-q}(G),J)
\qquad\text{and}\qquad
\DD 2pq=\Hom_R(\HH_{-q}(G),\dd^J_p)
 \]
and differentials acting by the pattern $\DD rpq\col\EE rpq\to
\EE r{p-r}{q+r-1}$ for all $r\ge2$.

If $J$ or $G$ is bounded below, then this sequence converges strongly
to $\Hom_R(G,J)$, in the sense that for each $n\in\BZ$ the group
$\HH_n\Hom_R(G,J)$ has a  finite, exhaustive, and separated filtration,
the component of degree $q$ of the associated graded group is isomorphic
to $\EE \infty {n-q}q$, and $\EE rpq=\EE \infty pq$ for all $r\gg0$
and all $(p,q)\in\BZ\times\BZ$.
 \end{proposition}

\begin{proof}
The filtration of $J$ by its subcomplexes $J_{\les p}$ induces a
filtration
 \[
\filt_p\Hom_R(G,J)=\Imag\big(\Hom_R(G, J_{\les p})\to \Hom_R(G,J)\big)
 \]
of $\Hom_R(G,J)$.  The differentials $\DD rpq$ in the resulting spectral
sequence $\EE rpq$ have the desired pattern for all $r\ge0$ and all
$(p,q)\in\BZ\times\BZ$.  The isomorphisms
 \[
\EE0pq={\Big(}\frac{\filt_p\Hom_R(G,J)}{\filt_{p-1}\Hom_R(G,J)}{\Big)}_{p+q}
\cong\Hom_R(G_{-q},J_p)
 \]
take $\DD 0pq$ to $(-1)^q\Hom_R(\dd^G_{-q+1},J_p)$.
As each $J_p$ is injective, the isomorphisms
 \[
\EE 1pq=\HH_{-q}\Hom_R(G,J_p)\cong\Hom_R(\HH_{-q}(G), J_p)
 \]
transfer $\DD 1pq$ into $\Hom_R(\HH_{-q}(G),\dd^J_p)$.
Thus, they lead to isomorphisms
 \[
\EE 2pq\cong\HH_p\Hom_R(\HH_{-q}(G),J)
 \]
that take $\DD 2pq$ to $\Hom_R(\HH_{-q}(G),\dd^J_p)$.

If $J$ or $G$ is bounded below, then the filtration of $\Hom_R(G,J)_n$ by
its subgroups $(\filt_p\Hom_R(G,J))_n$ is finite for each $n\in\BZ$; this
implies strong convergence.
 \end{proof}

Next comes a small variation on \cite[4.4]{AF}.

 \begin{lemma}
 \label{evaluation}
If $E$ is a complex of right $R$-modules, $G$ denotes the complex
of $R$-modules $\Hom_R(E,R)$, and $J$ is a complex of $R$-modules, then
the map
 \begin{align*}
\theta^{EJ}\col E\otimes_R J&\lra \Hom_R(G,J)
 \\
x\otimes y&\longmapsto\big(\gamma\mapsto
(-1)^{|\gamma||y|}\gamma(x)y\,\big)
 \end{align*}
is a morphism of complexes.  It is bijective if $J$ is bounded above,
$E$ is bounded below, and each right module $E_h$ is finite projective.
 \end{lemma}

 \begin{proof}
It is easily verified that $\theta^{EJ}$ is a morphism of complexes.

By definition, $\Hom_R(G,J)_n=\prod_{i-(-h)=n} \Hom_R(G_{-h}, J_i)$ for
each $n\in\BZ$.  The boundedness hypotheses imply that the abelian groups
$\Hom_R(G_{-h}, J_i)$ are trivial for almost all pairs $(h,i)$ satisfying
$i-(-h)=n$, so the product of these groups is equal to their coproduct.
The homomorphism $\theta^{EJ}_n\col(E\otimes_RJ)_n\to\Hom_R(G,J)_n$
is thus the coproduct of the canonical evaluation maps $E_h\otimes_RJ_i
\to\Hom_R(G_{-h},J_i)$, each one multiplied by $(-1)^{ni}$.  For every
$h\in\BZ$ the $R$-module $G_{-h}=\Hom_R(E_h,R)$ is finite projective,
all the evaluation homomorphisms are bijective.
 \end{proof}

Now we are ready to produce some canonical isomorphisms.  Our approach
was suggested by results in an early version of a paper of Huneke and
Leuschke \cite{HL}.

When $N$ is a right $R$-module $N^*=\Hom_R(N,R)$ carries the standard
left action$.$

 \begin{proposition}
 \label{ss}
Let $E$ be a complex of finite projective right $R$-modules with $\inf
E=0$, let $J$ be a complex of $R$-modules, and form the composition of
morphisms
 \[
\vartheta^{EJ}\col E\otimes_RJ\xra{\theta^{EJ}}\Hom_R(\Hom_R(E,R),J)
\xra{\Hom_R(\Hom_R(\epsilon,R),J)}\Hom_R(N^*,J)
 \]
where $N=\HH_0(E)$, the map $\epsilon\col E\to N$ is canonical, and
$\theta^{EJ}$ is from Lemma {\em\ref{evaluation}}.

When $J$ is bounded and each $R$-module $J_i$ is injective the
following hold.
 \begin{enumerate}[\quad\rm(1)]
  \item
If $\HH_{-i}\Hom_R(E,R)=0$ for all $i>0$, then $\vartheta^{EJ}$ is a
quasi-isomorphism.
  \item
If $\HH_{-i}\Hom_R(E,R)=0$ for all $i\gg 0$, then
$\HH_i(E\otimes_RJ)=0$ for all $i\gg 0$.
  \item
If $\sup J=0$ and there is number $m\ge 1$ such that $\HH_{-i}\Hom_R(E,R)=0$
for all $i\in[1,m]$, then for every $i\le m+\inf J$ there is an
isomorphism
 \[
\HH_i(\vartheta^{EJ})\col\HH_i(E\otimes_RJ)\xra{\ \cong\ }
\HH_i\Hom_R(N^*,J)
 \]
 \end{enumerate}
  \end{proposition}

 \begin{proof}
For $G=\Hom_R(E,R)$ Lemma \ref{evaluation} yields an
isomorphism of complexes
 \[
\theta^{EJ}\col E\otimes_RJ\xra{\ \cong\ } \Hom_R(G,J)
 \]

(1) By hypothesis, $\Hom_R(\epsilon,R)\col N^*\to G$
is a quasi-isomorphism, hence by \ref{semiinjective} so is
$\Hom_R(\Hom_R(\epsilon,R),J)$; as $\theta^{EJ}$ is bijective, we
are done.

(2) The hypothesis means that $\HH(G)$ is bounded.  The spectral sequence
of Proposition \ref{spectral sequence} shows that so is $\HH\Hom_R(G,J)$,
which is isomorphic to $\HH(E\otimes_RJ)$.

(3) Proposition \ref{spectral sequence} yields a strongly convergent
spectral sequence
 \[
\EE 2pq=\HH_p\Hom_R(\HH_{-q}(G),J)\implies\HH_{p+q}\Hom_R(G,J)
 \]
Note that $\inf J=-n$ for some integer $n\ge0$.  Because $G_{-q}=0$ for
all $q<0$, and $\HH_{-q}(G)=0$ for all $q\in[1,m]$, we get
$\EE 2pq=0$ for all pairs $(p,q)$ with $p\notin[-n,0]$ or with $q\le m$,
unless $q=0$.  It follows that if $p+q\le m-n$, then
 \[
\EE\infty pq=\EE 2pq\cong \begin{cases}
\HH_p\Hom_R(N^*,J) &\text{when } q=0\\
0 &\text{when } q\ne 0
 \end{cases}
 \]
As a consequence, the edge map $\HH_p\Hom_R(G,J)\to\EE 2p0$ is bijective
for $p\le m-n$.  It follows from the construction of the spectral sequence in
Proposition \ref{spectral sequence} that the isomorphism above takes this
edge map to $\HH_p\Hom_R(\Hom_R(\epsilon,R),J)$.
 \end{proof}

The next corollary contains as special cases \cite[2.3]{HH} and part
of \cite[2.1]{HH}.

\begin{corollary}
 \label{vs}
Let $L$ be a right $R$-module that has a projective resolution $E$
by finite projective right $R$-modules, and let $M$ be an $R$-module
with $\id_RM=n<\infty$.

If $m\ge 1$ is a number such that $\Ext^i_R(L,R)=0$ for all $i\in[1,m]$,
then
 \begin{enumerate}[\quad\rm(1)]
\item $\Ext^i_R(L^*,M)=0$ for all $i>\max\{0,n-m-1\}$.
\item $\Tor_i^R(L,M)=0$ for all $i\in[1,m]$.
\item If $m\ge n$, then $L\otimes_R M\cong\Hom_R(L^*,M)$.
 \end{enumerate}
 \end{corollary}

 \begin{proof}
Let $M\to J$ be an injective resolution with $J_i=0$ for $i<-n$.
The last part of the proposition yields for all $i\le m-n$ isomorphisms
 \[
\Tor_i^R(L,M)\cong\Ext^{-i}_R(L^*,M)
 \]
The desired assertions follow,  since $\Tor_j^R(L,M)=0=\Ext^j_R(L^*,M)$
for $j<0$.
 \end{proof}

Here is a version of the familiar degree-shifting procedure.

 \begin{lemma}
 \label{shifts}
Let $L$ be a complex of right modules with $\sup\HH(L)=l<\infty$ and
$M$ a complex with $\sup\HH(M)=m<\infty$.  If $E\to L$ and $F\to M$ are
semiprojective resolutions, then for the $R$-modules $L'=\Coker(\dd^E_l)$
and $M'=\Coker(\dd^F_m)$ there are isomorphisms $\Tor^R_i(L,M)\cong
\Tor^R_{i-l-m}(L',M')$ for all $i>l+m$.
 \end{lemma}

 \begin{proof}
The quasi-isomorphism $L\simeq \tau_{\les l}(L)$, cf.\ \ref{smart},
induces by \eqref{semiflat} a quasi-isomorphism $L\otimes_RF_{<m} \simeq
\tau_{\les l}(L)\otimes_R(F_{<m})$.  As a consequence, we get
 \begin{align*}
\sup \HH(L\otimes_R(F_{<m}))&=\sup\HH(\tau_{\les l}(L)\otimes_R(F_{<m}))\\
&\le \sup (\tau_{\les l}(L)\otimes_R(F_{<m}))\\ &= l+m-1
 \end{align*}
The homology exact sequence of the exact sequence of complexes
 \[
0\to L\otimes_R(F_{< m})\lra L\otimes_RF\lra L\otimes_RF_{\ges m}\to 0
 \]
now yields isomorphisms $\HH_i(L\otimes_RF)\cong \HH_i(L\otimes_RF_{\ges
m})$ for all $i>l+m$.  Since $\shift^{-m}F$ is a projective
resolution of $M'$, these isomorphisms can be rewritten as
$\Tor^R_i(L,M)\cong\Tor^R_{i-m}(L,M')$ for all $i>l+m$.  Similar
arguments applied to the resolution $E\to L$ yield $\Tor^R_{i-m}(L,M')
\cong\Tor^R_{i-l-m}(L',M')$ for all $i>l+m$.
 \end{proof}

\end{appendix}

\section*{Acknowledgements}

The authors want to thank Craig Huneke and Graham Leuschke for making
an early version of \cite{HL} available to them, Bernd Ulrich for
useful discussion concerning the material of Sections \ref{Well}
and \ref{Residue}, and the referee for directing them to Proposition
\ref{associated}.

 \end{document}